\documentclass[a4paper, 12pt]{article}
\usepackage[utf8]{inputenc}
\usepackage{makeidx}
\usepackage{enumerate, theorem}
\usepackage{amsmath, amsfonts, amssymb}
\usepackage[height=22.5cm, width=15cm]{geometry}
\usepackage[all]{xy}
\usepackage{mathrsfs, yfonts}
\usepackage{graphicx}
\DeclareMathOperator{\C}{\mathbb{C}}
\newcommand{\parag}[1]{\paragraph{\sc{#1.}}}
\newcommand{\Cech}{\check}

\newtheorem{thm}{Theorem}[subsection]

\newtheorem{cor}[thm]{Corollary}
\newtheorem{prop}[thm]{Proposition}
\newtheorem{lemma}[thm]{Lemma}

\begin{document}

\title{On Symmetric and Anti-symmetric Partial Differential Operators}
\date{11/11/2024}
 \author{Daniel Barlet\footnote{Institut Elie Cartan, G\'eom\`{e}trie,\newline
Universit\'e de Lorraine, CNRS UMR 7502   and  Institut Universitaire de France.}.}

\maketitle

\begin{align*}
&\qquad \qquad \qquad\qquad \qquad \qquad\qquad \qquad \qquad  & {\it Sym{\acute e}trie, \ vous  \ avez \  dit \ sym{\acute e}trie \ ??} \\
&\qquad \qquad \qquad\qquad \qquad \qquad\qquad \qquad \qquad  & {\it  Alors \ ne \ faites \ pas \ aux \ autres \ ce  \ que} \\
& \qquad \qquad \qquad\qquad \qquad \qquad\qquad \qquad \qquad  & {\it vous  \ ne \ voulez \  pas \ que \ l'on \ vous \ fasse.}
\end{align*}

\bigskip

\parag{Abstract} We study the action of symmetric PDO on the module of anti-symmetric functions in $z_1, \dots, z_k$. We show that over the Weyl algebra of the elementary symmetric functions, this module is simple.

\parag{AMS Classification} 32 C 38; 14 F 10

\parag{Key words} Symmetric partial differential operators- Antisymmetric polynomials-Discriminant.

\tableofcontents

\section{Introduction} 

This paper studies the module of antisymmetric functions $\C[\sigma]\delta$ in the variables  $z_1, \dots, z_k$ (where $\delta$ is the discriminant) under the action of the corresponding  Weyl algebra  $W_1^{\frak{S}_k} := \C \langle z_1, \dots, z_k \rangle^{\frak{S}_k}$ of symmetric partial differential operators. We also study   the action  of the Weyl algebra $W_2 := \C\langle \sigma_1, \dots, \sigma_k \rangle$ associated to $\sigma_1, \dots, \sigma_k$  the elementary symmetric functions of $z_1, \dots, z_k$, on $\delta$ inside  $\C[\sigma, \Delta^{-1}]\delta$.\\
Our first result (Theorem \ref{generators 1}) is to show that the algebra $W_1^{\frak{S}_k}$ is generated by $N_1, N_k$ and $\mathcal{N}_k$ which are respectively the Newton symmetric functions of $z_1, \dots, z_k $ and of $\partial/\partial z_1, \dots, \partial/\partial z_k$. This allows to determine the left ideal $\mathcal{I}$  in $W_1^{\frak{S}_k}$ which is the annihilator of $\delta$.  Then we show (Theorem \ref{simple}) that $\mathcal{I}$ is a maximal left ideal in $W_1^{\frak{S}_k}$, so the $W_1^{\frak{S}_k}$-module $\C[\sigma]\delta$ is simple.  As a corollary, we obtain that the left ideal $\mathcal{J}$ in $W_2$ which is the annihilator of $\delta$ is also a maximal left ideal in $W_2$.\\
So the $W_2$-module $W_2\delta := W_2/\mathcal{J} \subset \C[\sigma][\Delta^{-1}]\delta$ is simple.\\
In the last section, we give the  description of  the symmetric and anti-symmetric vector fields in $W_1$ and we show that any anti-symmetric differential operator $A$  belongs to $\delta W_2$. \\
As an application we show that this implies that the Bernstein polynomial of  the polynomial $\Delta = \delta^2$ in $\C[\sigma]$  does not vanish at $\lambda = -1/2$.

\bigskip

\section{Symmetric differential operators}

\subsection{Generators of $W_1^{\frak{S}_k}$}

\parag{Notations} Let $z_1, \dots, z_k$ be $k$ complex variables and let $\partial/\partial z_j, j \in [1,k]$, denote the corresponding partial derivatives.\\
Let $\sigma_1,\dots, \sigma_k$ the elementary symmetric functions of $z_1, \dots, z_k$ and $N_p, p \geq 0$, the corresponding Newton  symmetric functions.\\
Then $\C[\sigma]$ is the $\C$-algebra of $\frak{S}_k$-invariant polynomials in $z_1, \dots, z_k$.\\
Let $\Sigma_1, \dots, \Sigma_k$ the elementary symmetric functions of $\partial/\partial z_1, \dots, \partial/\partial z_k$ and let  $\mathcal{N}_p, p \geq 0$, be  the corresponding Newton symmetric functions.\\
We denote by $\C[\Sigma]$ the (commutative) algebra generated by $\Sigma_1, \dots, \Sigma_k$ which contains the $\mathcal{N}_p,\ \forall  p \geq 0$. \\
We denote respectively by $W_1$ and $W_2$ the Weyl algebras $\C \langle z_1, \dots, z_k \rangle$ and $\C \langle \sigma_1, \dots, \sigma_k \rangle$. There is a natural morphism of algebras $q_* : W_1 \to W_2$ where  we denote by $q$  the quotient  map  $q : M = \C^k \to N = \C^k$
by the natural  action of the permutation group $\frak{S}_k$ on $\C^k$. The restriction of $q_*$ to the sub-algebra $W_1^{\frak{S}_k}$ of symmetric  elements in $W_1$ is injective.\\
We denote by $\delta := \prod_{1\leq i < j \leq k} (z_i - z_j)$  and by $\Delta := \delta^2 \in \C[\sigma]$ the discriminants.\\

Our first result is to give a set of generators for the unitary $\C$-algebra $ W_1^{\frak{S}_k}$ of symmetric partial differential operators inside the Weyl algebra $W_1 := \C\langle z_1, \dots, z_k\rangle$.

\begin{thm}\label{generators 1}
The unitary $\C$-algebra $W_1^{\frak{S}_k}$ is generated  by $N_1, N_k$ and $\mathcal{N}_k$.
\end{thm}

The proof will use the following proposition where we define the algebra $\mathcal{A}$ as the sub-algebra in $W_1^{\frak{S}_k}$ generated by $\sigma_1, \dots, \sigma_k$ and $\Sigma_1, \dots, \Sigma_k$ which are respectively the elementary symmetric functions of $z_1, \dots, z_k$ and $\partial/\partial z_1, \dots, \partial/\partial z_k$.

\begin{prop}\label{23/10}
Let $p$ and $q$ be non negative integers. Define
$$ V_{p,q} := \sum_{j=1}^k  z_j^p(\frac{\partial}{\partial z_j})^q .$$
Then the $V_{p, q}$ generates the algebra $\mathcal{A}$.
\end{prop}

\parag{Proof} First remark that $V_{p, 0} = N_p$ and $V_{0, q} = \mathcal{N}_q$ are respectively the $p$-th and $q$-th Newton symmetric function of $\sigma_1, \dots, \sigma_k$ and $\Sigma_1, \dots, \Sigma_k$, so belongs to $\mathcal{A}$. This already shows that $\mathcal{A}$ is contained in the sub-algebra generated by the $V_{p, q}$.\\
In order to compute the commutators of two differential operators  $V_{p,q}$, remark that $z_j^p(\partial/\partial z_j)^q$ and $z_i^{p'}(\partial/\partial z_i)^{q'}$ commute for $i \not= j$; so the main point is to compute the commutator
of $[\partial^q, z^p] = pq z^{p-1}\partial^{p-1} + \cdots$. This may help the reader  to check the formulas we use below.\\
To prove the converse, that is to say that each  $V_{p, q}$ is in $\mathcal{A}$, remark that  the relation $4V_{1, 1} = [V_{2, 0}, V_{0, 2}] - 2V_{0, 0}/k$ \  gives $V_{1, 1} \in \mathcal{A}$.  Also the relations:
\begin{align*}
&  [V_{0, 3}, V_{2, 0}] = 6V_{1, 2}  +  6V_{0, 1}, \quad [V_{0, 2}, V_{3, 0}] = 6V_{2, 1} + 6V_{1, 0},\\
& \quad {\rm and} \quad  [V_{1, 2}, V_{2, 1}] = 3V_{2, 2} +2V_{1, 1}, 
\end{align*}
implies that $V_{p, q}$ is in $\mathcal{A}$ for any $p \leq 2$ and $q \leq 2$.\\
Now let $a \geq 2$ and assume that we have proved that $V_{p, q}$ is in $\mathcal{A}$ for any $p \leq a$ and $q \leq a$. Then the relations for $p \leq a$ and any  $q \leq a+1$.
\begin{align*}
& [V_{1, 2}, V_{p, a}] - (2p-a)V_{p, a+1} \in \mathcal{A} \quad {\rm and} \\
& [V_{2, 2}, V_{p, a}] - 2(p - a)V_{p+1, a+1} \in \mathcal{A} \quad  {\rm to \ treat \ also \ the \ case} \ V_{a/2, a+1} \\
& [V_{2, 1}, V_{a, q}] - (a - 2q)V_{a+1, q} \in \mathcal{A}  \quad {\rm and} \\
& [V_{2, 2}, V_{a, q}] - 2(a - q)V_{a+1, q} \in \mathcal{A} \quad  {\rm to \ treat \ also \ the \ case} \ V_{a+1, a/2} 
\end{align*}
shows that $V_{p, q}$ is also in $\mathcal{A}$ for each $p $ and $q$ at most equal $a+1$ completing our induction and the proof.$\hfill \blacksquare$\\

\parag{Remark} Consider the Weyl algebra $W := \C \langle z, \partial/\partial z \rangle$ in one variable $z$ with basis $e_{p, q} := z^p(\partial/\partial z)^q$ and define the $\C$-linear isomorphism
$\varphi : W \to W_1^{\frak{S}_k}$ sending $e_{p, q}$ onto $V_{p, q}$. Then this induces an isomorphism of (bi-graded) Lie  algebras. 
This an obvious consequence of the  fact that $z_j^p(\partial/\partial z_j)^q$ and $z_i^{p'}(\partial/\partial z_i)^{q'}$ commute for $i \not= j$:
$$ [\varphi(e_{p, q}), \varphi(e_{p', q'})] = [V_{p, q}, V_{p', q'}] = \sum_{j=1}^k [z_j^p(\partial/\partial z_j)^q, z_j^{p'}(\partial/\partial z_j)^{q'}] = \varphi ([e_{p, q}, e_{p', q'}]). $$
But $\varphi$ is not an algebra morphism since, for instance, $ N_pN_q \not= N_{p+q}$ !\\
Using the description of symmetric vector fields in $W_1$ (recalled in Section 3),  we see that $\varphi$ induces a Lie algebra isomorphism between polynomial vector fields on $\C$ and symmetric polynomial vector fields in $W_1$.$\hfill \square$\\

 We shall prove Theorem \ref{generators 1} by showing that   the sub-algebra $\mathcal{A}$ generated by the $V_{p, q}$  is equal to $W_1^{\frak{S}_k}$. So the conclusion follows thanks to Proposition \ref{23/10}

\parag{Proof of the theorem} Remark first that $\Sigma_1 = \mathcal{N}_1$ is equal to $V_{0,1} = [\mathcal{N}_2, N_1]$. \\
We shall prove that for any $(\alpha, \beta) \in (\mathbb{N}^k)^2$ the symmetrization of the monomial differential operator $z^\alpha (\partial /\partial z)^\beta$ is in the sub-algebra $\mathcal{A}$.\\
Assume that $h$ is an integer in $[0, k]$ and that the monomial differential operator $z^\alpha (\partial /\partial z)^\beta$ depends only on the variables $z_1, \dots, z_h$.\\

Then for $h = 0$ its symmetrization is in $\C[\sigma] \subset \mathcal{A}$ and for $h = 1$  its symmetrization is $(k-1)!V_{\alpha_1,\beta_1}$ which is also in $\mathcal{A}$ thanks to 
Proposition \ref{23/10}. \\
So assume that $h$ is in $[1,k-1]$ and that it is already proved that the symmetrization
 $$\frak{S}(z^\alpha (\partial /\partial z)^\beta) := \sum_{\tau \in \frak{S}_k} z^{\tau(\alpha)} (\partial /\partial z)^{\tau(\beta)} $$
 of $z^\alpha (\partial /\partial z)^\beta$  is in $\mathcal{A}$ for any $(\alpha,\beta) \in (\mathbb{N}^h)^2$. \\
  Then consider  $(\alpha,\beta) \in (\mathbb{N}^{h+1})^2$ and write
 $$ z_1^{\alpha_1}\cdots z_{h+1}^{\alpha_{h+1}}(\partial /\partial z_1)^{\beta_1}\cdots (\partial /\partial z_{h+1})^{\beta_{h+1}} = \big(z^{\alpha'}(\partial /\partial z)^{\beta'} \big) z_{h+1}^{\alpha_{h+1}}(\partial /\partial z_{h+1})^{\beta_{h+1}} $$
 where $(\alpha',\beta')$ is in $(\mathbb{N}^h)^2$, using the fact that $z_{h+1}$ and $\partial /\partial z_{h+1}$ commutes with $z_1, \dots, z_h$ and $\partial /\partial z_1, \dots, \partial /\partial z_h$.\\

Then product
$$ \frak{S}(z^{\alpha'}(\partial /\partial z)^{\beta'})V_{\alpha_{h+1},\beta_{h+1}} $$
 is in  the sub-algebra $\mathcal{A}$ thanks to our inductive hypothesis and to Proposition \ref{23/10}. This element in $\mathcal{A}$ is sum of monomials like
$$ z^{\tau(\alpha')}(\partial /\partial z)^{\tau(\beta')})z_j^{\alpha_{h+1}}(\partial /\partial z_j)^{\beta_{h+1}} $$
with $\tau \in \frak{S}_k$ and $j \in [1,k]$. Then there are two cases.
\begin{enumerate}
\item The case where $j$ is in $\{\tau(1), \dots, \tau(h)\}$.
\item The case where $j$ is not in $\{\tau(1), \dots, \tau(h)\}$ and then the corresponding monomial appears in $\frak{S}(z^\alpha (\partial /\partial z)^\beta) $.
\end{enumerate}
The sum of  the monomials in the first case is a symmetric partial differential operator which is the symmetrization of monomials involving at most $h$ variables, so is in $\mathcal{A}$ by our inductive assumption.\\
Looking at the coefficient of $z^{\tau(\alpha)}(\partial /\partial z)^{\tau(\beta)}$ in this sum we find that it appears exactly  $d$ times, if $d$ is the coefficient of  $z^{\tau(\alpha')}(\partial /\partial z)^{\tau(\beta')}$
in $\frak{S}(z^{\tau(\alpha')}(\partial /\partial z)^{\tau(\beta'})$ since we must choose $j = \tau(h+1)$. Note that $d$ is independent of $\tau \in \frak{S}_k$ and depends only on $(\alpha', \beta')$. If $c$  is the coefficient of $z^{\tau(\alpha)}(\partial /\partial z)^{\tau(\beta)}$  in $\frak{S}(z^\alpha (\partial /\partial z)^\beta) $. Also $c$ depends only on $(\alpha, \beta)$ and not on the choice of $\tau$. So we obtain that the sum of terms in case 2 is equal to 
$$\frac{d}{c}\frak{S}\big(z^\alpha(\partial /\partial z)^\beta\big).$$
So $\frak{S}(z^\alpha (\partial /\partial z)^\beta)$ is in $\mathcal{A}$ and  this completes our induction step.\\
To complete the proof, since we know now that $N_0, N_1, \dots, N_k$ and $\mathcal{N}_1, \dots, \mathcal{N}_k$ is a generator of the algebra  $\mathcal{A} = W_1^{\frak{S}_k}$, it is enough to use the equalities:
$$ [\mathcal{N}_h, N_1] = h\mathcal{N}_{h-1} \quad {\rm and} \quad [\mathcal{N}_1, N_h] = hN_{h-1} $$
which allows to obtain that $\mathcal{N}_1, \dots, \mathcal{N}_{k-1}$ and   $N_0, N_1, \dots, N_{k-1}$  are in the algebra generated by $N_1, N_k$ and $\mathcal{N}_k$.$\hfill\blacksquare$\\

\subsection{Generators for $Ann(\delta)$}

\parag{The map $\varphi$} Let $P \in W_1^{\frak{S}_k}$ and let $\delta := \prod_{1 \leq i < j \leq k} (z_i - z_j)$ be the discriminant. The polynomial $P(\delta) \in \C[z_1, \dots, z_k]$ is anti-symmetric so it may be written
in an unique way as $P(\delta) = \varphi(P).\delta$ where $\varphi(P) $ belongs to $ \C[z_1, \dots, z_k]^{\frak{S}_k} \simeq \C[\sigma_1, \dots, \sigma_k] $, since the $\C[\sigma]$-module of anti-symmetric polynomials in 
$\C[z_1, \dots, z_k]$ is equal to $\C[\sigma].\delta$. This defines a left  $\C[\sigma]$-linear map
$$ \varphi : W_1^{\frak{S}_k} \to \C[\sigma] .$$
This map has the following properties
\begin{enumerate}
\item $\varphi(PQ) = \varphi(P.\varphi(Q))$ for any $P$ and $Q$ in $W_1^{\frak{S}_k}$.
\item The annulator of $\delta$ in $W_1^{\frak{S}_k}$ is $Ker \varphi$ which is the left ideal is given by 
$$ Ker \varphi := \{ P -  \varphi(P) \ / \ P \in W_1^{\frak{S}_k}\}.$$
\item Let  $w$ be the weight on $W_1^{\frak{S}_k}$ defined by $w(z_j) = 1$ and $w(\partial/\partial z_j) = -1$ for each $j \in [1, k]$. The map $\varphi$ preserves pure weights. This is clear because $\delta$ has pure weight.
\item As a consequence \ $Ker \varphi$ \ contains any $P$ with pure negative weight. \\
For instance any  element without constant  in $\C[\Sigma]$ is in $Ker \varphi$.
\item If $V$ is a vector field and $P$ any element in  $W_1^{\frak{S}_k}$ we have
 $$ \varphi(VP) = V(\varphi(P)) + \varphi(P)\varphi(V).$$
So for two vector fields $V$ and $W$ we obtain $\varphi([V, W]) = V(\varphi(W)) - W(\varphi(V))$.
\item Since the vector field $V_{0, 1} = \Sigma_1 = \mathcal{N}_1 $ which is of weight $-1$  kills $\delta$ we have  $\varphi(\Sigma_1 P) = \Sigma_1(\varphi(P))$ for each $P \in W_1^{\frak{S}_k}$. 
\item For each $P \in W_1^{\frak{S}_k}$ we have $P(\Delta) = \varphi(\delta P\delta)$ where $\Delta := \delta^2$:\\
Since $\delta P \delta$ is symmetric for $P$ symmetric, we have
 $$(\delta P \delta)(\delta) = \delta P(\Delta) = \varphi(\delta P \delta)\delta$$
 and the conclusion follows.
\end{enumerate}

\begin{prop}\label{Annul.1}
For each $(p,q) \in \mathbb{N}$ there exists a polynomial $u_{p,q} \in \C[\sigma]$ of weight $p - q$  such that 
$$V_{p,q}(\delta) = u_{p,q}\delta \quad {\rm so} \quad \varphi(V_{p, q}) = u_{p, q}.$$
 Moreover we have
\begin{enumerate}[(i)]
\item For \  $p - q < 0$ \ or  for  \ $q > k, \ u_{p,q} = 0$.
\item For $q = 0, \  u_{p, 0} = N_p$.
\item For $q = 1, \ 2 u_{p+1, 1} = \sum_{h=0}^p N_hN_{p-h} - (p+1)N_p.$
\item For any $(p,q)$ the following formulas holds, for $R \gg \vert\vert \sigma \vert\vert $ where we put \\ $ \vert\vert \sigma \vert\vert  := Sup\{ \vert \sigma_h\vert^{1/h},   h \in [1, k]\}$:
\begin{align}
& (q+1)u_{p,q} = \frac{1}{2i\pi}\int_{\vert \zeta\vert = R} \ \frac{\zeta^p P^{[q+1]}(\zeta)}{P(\zeta)}d\zeta  \\
& (q+1) u_{p,q}  = \sum_{j=1}^k \frac{z_j^pP^{[q+1]}(z_j)}{P'(z_j)} \\
&  (q+1)u_{p,q} = \sum_{h=0}^{p-q} (-1)^h \frac{(k-h)!}{(k-h-q-1)!}\sigma_hM_{p-q-h}
 \end{align}
 where $P^{[h]}$ is the $h$-th derivative of $P(z) := \prod_{j=1}^k (z-z_j)$, and where we define for $i \geq 1 - k$:
 $$M_i := \sum_{j=1}^k \frac{z_j^{k+i-1}}{P'(z_j)} = (-1)^{h-1} \frac{1}{i+h}\partial_h N_{i+h}$$
 for each $h \in [1,k]$,  and where $\partial _1, \dots, \partial_k$ are the partial derivatives in the coordinates  $\sigma_1, \dots, \sigma_k$.
 \end{enumerate}
 \end{prop}
 
 The proof of this proposition uses the following lemma.
 
  \begin{lemma}\label{correction}
 Let $P(z) := \prod_{h=1}^k (z- z_h) $ and for $j \in [1, k]$ put $$\Pi_j(z) := \prod_{h\not= j, h=1}^k (z - z_h).$$
  Then we have $P'(z_j) = \Pi_j(z_j)$ and for each $q \in \mathbb{N}$ 
 \begin{equation*}
  \partial^q(\Pi_j(z_j))/\partial^q z_j = \frac{1}{q+1}P^{[q+1]}(z_j). \tag{$@$}
  \end{equation*}
 \end{lemma}
 
 Note that $\Pi_j$ is a polynomial in $z$ which is independent of $z_j$. So $\Pi'_j(z_j) = (\partial /\partial z_j)(\Pi_j(z_j))$.
 
 \parag{Proof} First note that $P'(z_j) = \Pi_j(z_j)$ for each $j \in [1, k]$. Then, let $ \sigma_h(\hat{j})$ be  the $h$-th symmetric function of $z_1, \dots, \hat{z_j}, \dots z_k$ and use the fact that 
 $$\sigma_h =  \sigma_h(\hat{j}) + \sigma_{h-1}(\hat{j})z_j, \ {\rm and \ so} \quad \frac{\partial \sigma_h}{\partial z_j} =  \sigma_{h-1}(\hat{j}) \quad  {\rm with} \ \sigma_0 (\hat{j}) \equiv 1$$
   to obtain
 \begin{align*}
 & \Pi_j'(z_j) = (\partial/\partial z_j)\big[\sum_{h=0}^{k-1} (-1)^h (k-h)\sigma_h z_j^{k-h-1}\big ]= P''(z_j) + \sum_{h=1}^{k-1} (-1)^h(k-h)\sigma_{h-1}(\hat{j})z_j^{k-h-1}\\
 &  \Pi_j'(z_j)  = P''(z_j) -\Pi_j'(z_j) \quad {\rm since } \quad  \Pi_j'(z) = \sum_{p=0}^{k-2} (-1)^p(k-p-1)\sigma_p(\hat{j})z^{k-p-2} .
 \end{align*}
 This gives the case $q = 1$ fo Formula $(@)$.\\
 Assume that we have proved that
 $$ (\partial^q/\partial z_j^q)[P'(z_j)] = \frac{1}{q+1}P^{[q+1]}(z_j) $$
 for some $q \geq1$. Then we obtain by derivation in $z_j$
 \begin{align*}
 & (\partial^{q+1}/\partial z_j^{q+1})[P'(z_j)] = \frac{1}{q+1}\Big(P^{[q+2]}(z_j) + \sum_{h=1}^{k-q-1} (-1)^h \frac{(k-h)!}{(k-h-q-1)!}\sigma_{h-1}(\hat{j})z_j^{k-h-q-1}\Big) \\
 & \qquad =  \frac{1}{q+1}\Big(P^{[q+2]}(z_j) - (\partial^{q+1}/\partial z_j^{q+1})[P'(z_j)] \Big) \quad {\rm since} \quad  P'(z_j) =  \Pi_j(z_j).
 \end{align*}
Now the equality $\big(1 + 1/(q+1)\big)(q+1) = q+2$ allows to conclude.$\hfill \blacksquare$\\

  \parag{Proof of Proposition \ref{Annul.1}} For any choice of $j \in [1, k]$, we may write 
   $$\delta = \prod_{1 \leq i < h \leq k} (z_i - z_h) = (-1)^{j-1}P'(z_j)\vartheta_j$$
   where $\vartheta_j$ is independent of $z_j$.  Then for any $q \geq 0$ we have, thanks to Lemma \ref{correction}
 $$ (\partial/\partial z_j)^q(\delta) = \frac{1}{q+1} \frac{P^{[q+1]}(z_j)}{P'(z_j)}\delta $$
 and then
 $$V_{p,q}[\delta] =  \frac{1}{q+1} \sum_{j=1}^k   \frac{z_j^p P^{[q+1]}(z_j)}{P'(z_j)}\delta \quad \forall (p, q) \in \mathbb{N}^2 .$$
The formulas of the proposition follow easily.$\hfill \blacksquare$\\
 
 The next result uses the same kind of technic that in  the proof of Theorem \ref{generators 1}
 
 \begin{prop}\label{Annul. 2}
 The left ideal in the algebra $W_1^{\frak{S}_k}$ of polynomial  symmetric differential operators on $\C^k$ which annihilates $\delta$ is generated by
 \begin{equation}
 V_{p,q} - u_{p,q} \quad {\rm for} \quad p \in [0,k], q \in  [1, k] 
 \end{equation}
 \end{prop}
 
 \parag{Proof} Remark first that it is enough to consider only pure weight element, since $\delta$ has pure weight $k(k-1)/2$. If $P$ has order $0$, then $P(\delta) = 0$ implies that $P = 0$. Then fix a weight $w$. We shall prove the proposition when $P$ has pure weight $w$. We shall argue by contradiction. For a non zero symmetric differential operator $P$ let $q$ its order, $h$ the maximal number of the variables $z_1, \dots, z_k$ which appears in a monomial $z^\alpha(\partial_z)^\beta$ (counting also differentials ) with a non zero coefficient in the order $q$ part of $P$ and let $\theta$ the number of such monomials in $z_1, \dots z_h$ (so we consider $P$ as the symmetrization of a sum of  such monomials). Then consider the lexicographical order on the set of triples $(q, h, \vartheta)$ and let  $(q_0, h_0, \theta_0)$ be minimal triple  for some $P$ which annihilates $\delta$ and which is not in the left ideal generated by the $V_{p,q} - u_{p,q}$ for all $(p,q) \in \mathbb{N}^2$ assuming that there exists such a $P$. Denote $P_0$ such an element giving this minimum.\\
 First remark that $q_0 \geq 1$ because no non zero element in $\C[\sigma]$ can annihilate $\delta$. Moreover, by definition $ h_0 \geq 1$ and $\theta_0 \geq 1$.\\
 Then consider one of the monomial  $z^\alpha(\partial_z)^\beta$  with $\vert \beta\vert = q$, involving only $z_1, \dots, z_{h_0}$ (we may assume that $h_0 \geq 2$ because when $h_0 = 1$ we may substract the corresponding $V_{\alpha_1, \beta_1} -u_{\alpha_1, \beta_1}$ and this monomial disappears !) and assume that $\beta_{1} \geq 1$. Then write
 $$ z^{\alpha}(\partial_z)^\beta = z^{\alpha '}(\partial_z)^{\beta '}z_1^{\alpha_1}(\partial_{z_1})^{\beta_1}$$ 
 where $\alpha'$ and $\beta'$ involve only $z_2, \dots z_{h_0}$. We obtain, $\frak{S}$ denoting the symmetrization operator:
 $$ \frak{S}(z^{\alpha}(\partial_z)^\beta) = \frak{S}(z^{\alpha '}(\partial_z)^{\beta '})(V_{\alpha_1, \beta_1} - u_{\alpha_1, \beta_1}) + \frak{S}(z^{\alpha '}(\partial_z)^{\beta '})u_{\alpha_1, \beta_1} + Q $$
 where $Q$ has numbers $(q, h, \theta)$ strictly less than $(q_0, h_0, \theta_0)$ and where $\frak{S}(z^{\alpha '}(\partial_z)^{\beta '})u_{\alpha_1, \beta_1} $ has order at most $q-1$.\\
 Indeed, any order $q_0$ monomial in $Q$ has at most $h_0$ different variables. Note that, if $h_0 = 1$   then $\beta ' = 0$ and  $q_0 = 1$,  case which is elementary because  we have\footnote{It is recalled in Section 3.1  that any symmetric vector field is the $\C[\sigma]$-module generated by the $V_{p,1}$ for $p \in [0, k-1]$.} 
 $ P = \big(\sum_{j=0}^w \mu_j V_{w+1-j,1}\big) + f_w$  where  the $\mu_j$ and   $f_w$ are in $\C[\sigma]$.\\
  So $P(\delta) = \big(f_w + \sum_{j=0}^w(\mu_j u_{w+1-j, 1} )\big)\delta = 0$ and then $P = \sum_{j=0}^w\mu_j(V_{w+1-j, 1} - u_{w+1-j, 1})$.\\
 Then
  $$P_1 := P_0 - \frak{S}(z^{\alpha '}(\partial_z)^{\beta '})(V_{\alpha_1, \beta_1} - u_{\alpha_1, \beta_1}) = \frak{S}(z^{\alpha '}(\partial_z)^{\beta '})u_{\alpha_1, \beta_1} + Q_0 $$
  has numbers $(q_1, h_1, \theta_1)$ with the following cases:
 \begin{enumerate}
 \item if $\theta_0 \geq 2$ then $\theta_1 = \theta_0 -1$,  and $h_1 = h_0, q_1 = q_0$.
 \item if $\theta_0 = 1$ then either $h_1 = h_0 - 1$ if $h_0 \geq 2$ or $q_1 = q_0 -1$ if $h_0 = 1$.
 \end{enumerate}
 In all case, we have  $P_1(\delta) = 0$ and $(q_1, h_1, \theta_1) < (q_0, h_0, \theta_0)$. So $P_1 $ is  in the left ideal generated by the $V_{p,q} - u_{p,q}$ and then $P_0$ also, which gives the contradiction.$\hfill \blacksquare$

 \parag{Remark}  \begin{enumerate}
 \item For each $p \geq k$ and each $q \geq 1$  we have
  $$V_{p,q} - u_{p,q} = \sum_{h=1}^k (-1)^{h-1} \sigma_h (V_{p-h,q}- u_{p-h,q}) $$ 
  and $V_{p,0} = u_{p,0} = N_p$ for each $p \geq 0$ because we have from $(1)$, the equality, for each $p \geq k$
  $$ \sum_{h=1}^k (-1)^{h-1}\sigma_hu_{p-h, q} = u_{p, q}.$$
 \item The elementary symmetric functions $\Sigma_1, \dots, \Sigma_k$ of $\partial/\partial z_1, \dots, \partial/\partial z_k$ are in the left ideal generated by $\mathcal{N}_1, \dots, \mathcal{N}_k$ the Newton symmetric functions of $\partial/\partial z_1, \dots, \partial/\partial z_k$ and $\mathcal{N}_h = V_{0,h}, \forall h \geq 1$ (note that $u_{0,h} = 0$ for $h \geq 1$).
 \item For each $q \geq k$ we have $V_{p,q} = \sum_{h=1}^k (-1)^{h-1}V_{p,q-h}\Sigma_h $.
 Note that it is clear, thanks to the previous lemma, that for $q \geq k$ we have $u_{p,q} = 0$.
 \item The differential operator $\mathcal{N}_2 := \sum_{j=1}^k (\partial/\partial z_j)^2$ is not in the left ideal of $W_1$ generated by the $V_{p,1}- u_{p, 1}$ for $p \in [0,k-1]$ since its symbol, $\sum_{j=1}^k \eta_j^2$ is not in the ideal of $\C[z_1,\dots, z_k, \eta_1, \dots, \eta_k]$ generated by the $\sum_{j=1}^k z_j^{p+1}\eta_j$. Indeed, if $\sum_{j=1}^ka_p^j(z)\eta_j$ are the homogeneous degree $1$ part in $\eta$ such that 
 $$ \sum_{p=-1}^{k-2} \Big[(\sum_{j=1}^k a_p^j(z)\eta_j)(\sum_{h=1}^k z_h^{p+1}\eta_h) \Big]= \sum_{j=1}^k \eta_j^2 $$
 we would obtain that $a_{-1}^j(z) \equiv 1$ for each $j \in [1,k]$ and $a_p^j(z) \equiv 0$ for $p$ in $[0, k-2]$ and any $j \in [1,k]$. This implies $\sum_{j=1}^k \eta_j^2 = \big(\sum_{j=1}^k \eta_j\big)^2$ which is absurd for $k \geq 2$.\\
 \item  Note that for each $j \in [1,k]$  the order $1$ differential operator 
 $$2P'(z_j)\partial /\partial z_ j - P''(z_j)$$
   kills  $\delta$.\\
 \end{enumerate}

 \subsection{Minimal extensions}
 
 We have seen above that $W_1^{\frak{S}_k}$ acts on $\C[\sigma_1, \dots, \sigma_k]\delta$ by the natural action of $W_1$ on polynomials in $\C[z_1, \dots, z_k]$. The next result shows that the corresponding  left $W_1^{\frak{S}_k}$-module is simple.
 
\begin{thm}\label{simple} The left  $W_1^{\frak{S}_k}$-module $\mathcal{M} := W_1^{\frak{S}_k}\big/\mathcal{I}$ where $\mathcal{I}$ is the annihilator of $\delta$ in $W_1^{\frak{S}_k}$ is a maximal left ideal in 
$W_1^{\frak{S}_k}$.
\end{thm}

The proof will use the following proposition.

  \begin{prop}\label{import}
Let $k \geq 2$ be an integer and let $\delta$ be the  discriminant in $z_1, \dots, z_k$.\\
If  some  $f \in \C[z_1, \dots, z_k]$ satisfies  $ \Sigma_h[f \delta] = 0 \quad \forall h \in [1, k]$, then  $f$ is constant.
\end{prop}

\parag{Remarks} \begin{enumerate}
\item We may replace the condition $\Sigma_h[f\delta] = 0,  \forall h \in [1, k]$ by analog  the condition:
 $\mathcal{N}_h[f\delta] = 0 , \ \forall h \in [1, k]$, because  $(\mathcal{N}_h, h \in [1, k])$  is also a generator of  the ideal generated by $\Sigma_1, \dots, \Sigma_h$ in the commutative algebra $\C[\Sigma_1, \dots, \Sigma_k]$.\\
\item Note that $\Sigma_h[\delta] = 0$ for each $h \in [1, k]$ since this anti-symmetric polynomial has degree strictly less than $w(\delta) = k(k-1)/2$, so $f = 1 $ satisfies the hypothesis.$\hfill \square$
\end{enumerate}

\parag{Proof} Let us begin by the case $k =2$. Let $f \in \C[z_1, z_2]$ be such that 
$$(\partial_{z_1} + \partial_{z_2})[(z_1 - z_2)f] = 0 \quad {\rm and} \quad  \partial_{z_1}\partial_{z_2}[(z_1 - z_2)f] =  0.$$
Then we have, using the second equation and $\partial_{z_1}\partial_{z_2} = \partial_{z_2}\partial_{z_1}$:
\begin{align*}
& (z_1 - z_2)\partial_{z_2}(f) - f (z_1, z_2)= g(z_2) \quad {\rm and} \\
& (z_1 - z_2)\partial_{z_1}(f) + f (z_1, z_2)= h(z_1)
\end{align*}
where $g, h$ are in $\C[z_1, z_2]$. For $z_1 = z_2 = z$ this gives  $f(z, z) = h(z) = -g(z)$ and  so   $(z_1 - z_2) \big(\partial_{z_1}+ \partial_{z_2}\big)[f](z_1, z_2) = g(z_2) - g(z_1)  = 0$ \  since \  $(\partial_{z_1} + \partial_{z_2})[(z_1 - z_2)f] = 0$.\\
So $g$ is a constant and  the  relation $\partial_{z_2}[(z_1 - z_2)f] = g$  gives $$(z_1 - z_2)f(z_1, z_2) = g z_2 + \gamma(z_1).$$
 Now $z_1 = z_2$ gives $\gamma(z) = -g z$ and  $f(z_1, z_2) = -g$, concluding the proof for $k = 2$.\\

For the case $k \geq 3$ remark that we may assume that $f$ is homogeneous in $z_1, \dots, z_k$.\\
We shall make  proof by an induction on $k \geq 2$ on the  the fact that if a homogeneous polynomial $f$ of degree $d$ satisfies $\Sigma_h[f\delta_k] = 0$ for each $h \in [1, k]$ then
 $d = 0$.\\
 So  we fix $k \geq 2$ and we assume that for $k-1$ we have proved that
$$ \Sigma_h[g \delta_{k-1}] = 0 \quad \forall h \in [1, k-1]  \quad {\rm implies} \quad  g \ {\rm is \ constant} $$
for $g \in \C[z_1, \dots, z_{k-1}]$ homogeneous.\\
Assume that the  polynomial $f$  in $z_1, \dots, z_k$ is homogeneous of degree $d$  and satisfies $\Sigma_h(f\delta_k) = 0$ for each $h \in [1, k]$.\\
 Then write, with $z' := (z_1, \dots , z_{k-1})$, 
\begin{align*}
& f(z) = \sum_{j=0}^d f_j(z')z_k^j \quad {\rm where} \quad f_j \in \C[z'] \quad {\rm is \ homogeneous \ with \   degree} \quad d-j \\
& \delta_k(z) = (-1)^{k-1} \Pi_k(z_k)\delta_{k-1}(z') \quad {\rm where} \quad \Pi_k(z) := \prod_{j=1}^{k-1} (z - z_j) =  \sum_{h=0}^{k-1} (-1)^h\sigma_h(z')z_k^{k-h-1}.
\end{align*}
Then put $f(z)\Pi_k(z_k) := \sum_{j=0}^{k+d-1} v_j(z')z_k^{d+k-j-1} $ where for each $j$ in  $[0, k+d-1]$  we have $v_j(z') = \sum_{h=0}^{inf(d, k-1)} (-1)^h\sigma_h(z') f_{d-j+h}(z')$ is homogeneous of degree $j$ in $z'$.\\
Then the equality $\Sigma_k[f\delta_k] = 0$ gives, since $\Sigma_k(z) = \Sigma_{k-1}(z')\partial_{z_k}.$\\
 We use here the notations $\Sigma_h(z)$ and $\Sigma_h(z')$ to distinguish the $h$-th symmetric functions of $\partial/\partial z_1, \dots, \partial/\partial z_k$ and $\partial/\partial z_1, \dots, \partial/\partial z_{k-1}$ respectively.\\
 Then
$$ \Sigma_k[f\delta_k] = (-1)^{k-1}\Sigma_{k-1}(z')\big[\delta_{k-1}(z')\partial_{z_k}(fP(z_k))\big] $$
implies the vanishing  for each $j \in [0, k+d-2]$ of the  coefficient of $z_k^{d+k-j-2}$ in  the polynomial $\Sigma_k[f\delta_k]$ in $z_k$, which is given by:
$$(d+k-j-1) (-1)^{k-1}\Sigma_{k-1}(z')[v_j(z')\delta_{k-1}(z')] = 0 .$$
So, for $j \in [0, d+k-2]$ we have $\Sigma_{k-1}(z')[v_j(z')\delta_{k-1}(z')] = 0. $ \\
Then the relation $\Sigma_{h}[f\delta_k] = 0$ gives, in the same way, the vanishing of 
\begin{align*}
&  \big(\Sigma_h(z') + \Sigma_{h-1}(z')\partial_{z_k}\big)[ f\Pi_k(z_k)\delta_{k-1}(z')] \quad {\rm which \ implies} \\
& \Sigma_h(z')[v_j(z')\delta_{k-1}(z')] + \Sigma_{h-1}(z')[(d+k-j-2)v_{j+1}(z')\delta_{k-1}(z')] = 0 \quad \forall j \geq 0 .
\end{align*}
As we already know that $\Sigma_{k-1}(z')[v_j(z')\delta_{k-1}(z')]  = 0$ for $j \in [0, d+k-2]$ this gives  for $h = k-1$ that
$$ 
\Sigma_{k-2}(z')[v_j(z')\delta_{k-1}(z')] = 0 \qquad \forall j \leq d+k-3.$$
Continuing this way we obtain 
$$ \Sigma_h(z')[v_j(z')\delta_{k-1}(z')] = 0 \quad \forall j \leq d+h-1 \quad {\rm for \ each} \quad h \in [1, k-1]$$
 Our inductive assumption implies then  that $v_j$ is constant for $j \in [0, d]$ and so  $v_j = 0$ for each  $1 \leq j \leq d$ and $v_0$ is constant.\\
 Then $f\Pi_k(z_k) = v_0z_k^{d+k-1} + R$ where $R$ has degree at most $k-2$ in $z_k$. This implies that $f = Q$ where $Q$ is the quotient of the division of 
 $v_0z_k^{d+k-1}$ by $\Pi_k(z_k)$.
 
 \bigskip
 
 To complete the proof of the proposition we  need some more results.
 
 \begin{lemma}\label{25/10 1}
 Consider now variable $z_1, \dots, z_{k-1}$ with elementary symmetric functions $s_1, \dots, s_{k-1}$ and  define  
 $\Pi_k(z) := \prod_{h=1}^{k-1} (z - z_h) = \sum_{p=0}^{k-1} (-1)^p s_p z^{k-p-1} $. For $d \in \mathbb{N}$ write the division of $z^{d+k-1}$ by $\Pi_k(z)$ as follows:
 \begin{equation}
 z^{d+k-1} = Q_d(z)\Pi_k(z) + R_d(z) \quad deg_z (R_d) \leq k-2 .
 \end{equation}
 Then we have, for $\rho$ large enough compared to $s_1, \dots, s_{k-1}$ and $z$:
 \begin{align*}
 & R_d(z) = \frac{-1}{2i\pi}\int_{\vert \zeta\vert = \rho}  \zeta^{d+k-1}\frac{\Pi_k(z) - \Pi_k(\zeta)}{\Pi_k(\zeta)(\zeta - z)} d\zeta \quad {\rm and} \\
 &  Q_d(z) = \frac{1}{2i\pi}\int_{\vert \zeta\vert = \rho}  \frac{\zeta^{d+k-1}}{\Pi_k(\zeta)(\zeta - z)} d\zeta.
 \end{align*}
 Moreover, for $z = z_k$ we have $Q_d(z_k) = M_d(\sigma)$ where the polynomial $M_d \in \C[\sigma]$ with weight $d$ is defined by
 \begin{align*}
 & M_d(\sigma) = \sum_{j=1}^k \frac{z_j^{d+k-1}}{P'(z_j)} =  \frac{1}{2i\pi}\int_{\vert \zeta\vert = \rho} \frac{\zeta^{d+k-1}}{P(\zeta)}d\zeta 
 \end{align*}
 where $P(\zeta) = \prod_{j=1}^k (\zeta - z_j) = \sum_{h=0}^k (-1)^h\sigma_h\zeta^{k-h}$ and $\rho \gg \vert\vert \sigma \vert\vert$.
 \end{lemma}
 
 \parag{Proof} This lemma is a standard consequence of Residues' formula (see, for instance,  \cite{[B.21]}  Section 3.4 for some  details).$\hfill \blacksquare$\\
 
 \begin{cor}\label{25/10 2}
 For each $d \in \mathbb{N}$ the polynomial $M_d$ (its definition is recalled in the previous lemma) is a monic polynomial of degree $d$ in $\sigma_1$ 
 (with coefficients in $\C[\sigma_2, \dots, \sigma_k]$) and it satisfies $\Sigma_1(M_d) = (d+k-1)M_{d-1}$ for each $d \geq 1$.\\
  So,  for $d \geq 1$,  we have  $\Sigma_1(M_d) \not= 0$.
 \end{cor}
 
 \parag{Proof} Recall that for each $p \in \mathbb{N}$ we have  $\partial_h(N_p) = (-1°^{h-1}pM_{p-h}$ where $N_p$ in $\C[\sigma]$ is the $p$-th Newton polynomial and where $M_d = 0$ for $d \in [-k+1, -1]$ (see \cite{[B.22]}). Since we have $\Sigma_1= \sum_{h=0}^{k-1} (k-h)\sigma_h\partial_{h+1}$ (with $\sigma_0 \equiv 1$), the commutation relations
 $$ [\Sigma_1, \partial_h] = -(k-h)\partial_{h+1} \quad \forall h \in [1, k-1]  \quad {\rm and} \quad [\Sigma_1, \partial_k] = 0$$
 hold true. So we have
 \begin{align*}
 & \Sigma_1\partial_1(N_{d+1}) = \partial_1\Sigma_1(N_{d+1}) - (k-1)\partial_2(N_{d+1}) \\
 & \Sigma_1((d+1)M_d) = \partial_1((d+1)N_d) + (k-1)(d+1)M_d = (d+1)(d+k-1)M_d
 \end{align*}
 proving our first assertion.\\
 Since $M_0 = 1, M_1 = \sigma_1$ and $M_{d} = -\sum_{h=1}^k (-1)^h\sigma_hM_{d-h}$ holds true for each  $d \geq 1$ (recall that $M_d = 0$ for $d \in [-k+1, -1]$), the last assertions follow.$\hfill \blacksquare$\\
 
 \parag{End of proof of Proposition \ref{import}} We have proved, using the previous lemma, that $f = v_0M_d(\sigma)$. But since $\Sigma_1(f\delta) = 0$ by assumption and  since
 $\Sigma_1(\delta) = 0$, we conclude that $f$ satisfies $\Sigma(f) = v_0 \Sigma(M_d) = 0$. For $d \geq 1$ the previous Corollary gives that $\Sigma_1(M_d) \not= 0$, so either $d = 0$ or $v_0 = 0$. In both cases $f$ is a constant, completing the proof.$\hfill\blacksquare$\\

\parag{Proof of the Theorem \ref{simple}} Define the weight of an element in $W_1^{\frak{S}_k} $ as the maximal weight of its monomials, and define the weight of $z^\alpha(\partial_z)^\beta$ as $\vert \alpha\vert - \vert \beta\vert$. Then recall  that the $\C[\sigma]$-linear map $\varphi : W_1^{\frak{S}_k} \to \C[\sigma]$ defined by $P(\delta) = \varphi(P)\delta$  keeps  the pure  weights since the weights of $\frak{S}(z^\alpha(\partial_z)^\beta)[\delta] $ is equal to $ \vert \alpha\vert - \vert \beta\vert + k(k-1)/2$  which implies $w(\varphi(P)) = w(P)$ for each pure weight  $P \in W_1^{\frak{S}_k}$ such that  $\varphi(P) \not= 0$.\\
Let $\mathcal{I}$  be the kernel of  the map $\varphi$ and choose  $P$ which  does not belong to $\mathcal{I}$.  Then put $ \mathcal{J} :=  \mathcal{I} + W_1^{\frak{S_k}}P$ and  choose  $P_0 \in \mathcal{J} \setminus \mathcal{I}$  with the condition that  $\varphi(P_0)$  has the smallest possible  weight in the image by $\varphi$ of $\mathcal{J} \setminus \mathcal{I} $.\\ Now we have, for each integer $h \geq 1$ 
$$ (\mathcal{N}_hP_0)(\delta)  = \varphi(\mathcal{N}_hP_0)\delta = \mathcal{N}_h[\varphi(P_0)\delta]  $$
 This shows that $\varphi(\mathcal{N}_hP_0) $ which has strictly smaller weight than $\varphi(P_0)$ must  vanish. Then the non zero element $\varphi(P_0)$ which has the smallest weight  in $\varphi(\mathcal{J})$ satisfies $\mathcal{N}_h(\varphi(P_0)\delta) = 0 $ for each $h \geq 1$. Proposition \ref{import} gives that $\varphi(P_0)$ is a constant (which is not $0$ since $P_0$ is not in $\mathcal{I}$).\\
 So, $1$ is in $\mathcal{J} = \mathcal{I} + W_1^{\frak{S}_k}P$ and $\mathcal{J} = W_1^{\frak{S}_k}$ proving that  any non zero  left sub-module of $W_1^{\frak{S}_k}$ containing  strictly $\mathcal{I}$  is equal to $W_1^{\frak{S}_k}$,  concluding the proof.$\hfill \blacksquare$\\

  We shall now define the action of  the Weyl algebra $W_2 := \C\langle \sigma, \partial \rangle $ on $\C[\sigma][\Delta^{-1}]\delta$ where $\Delta = \delta^2$.\\
  
  For $i \in [1, k]$ we define
  $$\partial_i(\delta) := \frac{\partial_i (\Delta)}{2\Delta}\delta $$
  and we extend this $\C[\sigma][\Delta^{-1}]$-connection on the rank $1$ free $\C[\sigma][\Delta^{-1}]$-module $\mathcal{M}$ with basis $\delta$ to a left  $W_2$-module structure on $\mathcal{M}$. \\
  
  Note that the restriction of this action to $W_1^{\frak{S}_k}$ induces an action on $\C[\sigma]\delta$ which is given by $P \mapsto (f\delta \mapsto P[f\delta])$ which sends $f\delta$ to  an anti-symetric polynomial in $\C[z_1, \dots, z_k]$ and so which is equal to $\varphi(Pf)\delta$ for some $\varphi(Pf) \in \C[\sigma]$. The $W_2$-connection on $\C[\sigma][\Delta^{-1}]\delta$ induces a left  $\C[\sigma][\Delta^{-1}]$-linear action of the localized algebra
   $$[\Delta^{-1}] W_2:= \cup_{m \in \mathbb{N}}\  (1/\Delta^m)W_2 = \C[\sigma][\Delta^{-1}] \otimes_{\C[\sigma]}W_2 $$
 on   $\mathcal{M}$.  Denote by $\mathcal{M}_0$ the left $W_2$-module generated by $\delta$ inside $\mathcal{M}$.\\
    We may define also the left $W_2$-module structure on $\mathcal{M}_0$ as follows, using the action of $W_1^{\frak{S}_k}$ on $\C[\sigma]\delta$:\\
  For $Q \in W_2$ there exists $m \in \mathbb{N}$ such that $\Delta^mQ$ belongs to $W_1^{\frak{S}_k}$(See Lemma \ref{14/8} in Section 3.1). Then define $Q(\delta)$ by the formula
   $$ Q(\delta) := \Delta^{-m}(\Delta^mQ)[\delta] \in \C[\sigma][\Delta^{-1}]\delta. $$
  Since for $P  \in W_1^{\frak{S}_k}$ we have $P[\delta] = \varphi(P)\delta$ which is in $\C[\sigma]\delta$, it is easy to see that this definition does not depend on the choice of $m$ such that $\Delta^mQ$ belongs to  $W_1^{\frak{S}_k}$. \\
  But with this definition it is clear that for $Q_1$ and $Q_2$ in $W_2$  the action of  $Q_2Q_1$ on $\delta$ is given by the action of $Q_2$ on
   $Q_1(\delta) \in \C[\sigma][\Delta^{-1}]\delta$ using $m$ large enough.
   
   \parag{Notation} We denote $[\Delta^{-1}] \mathcal{I} := \cup_{m \geq 0} \Delta^{-m}\mathcal{I}  \subset [\Delta^{-1}]W_2$.
  
  \begin{thm}\label{simple aussi} Let $\mathcal{J}$ be the left ideal in $W_2$ which is  the annihilator of $\delta$  in $\C[\sigma][\Delta^{-1}]\delta$. Then $\mathcal{J} = [\Delta^{-1}] \mathcal{I} \cap W_2$ and $\mathcal{J}$ is a maximal left ideal in $W_2$. So the left $W_2$-module $\mathcal{M}_0$ is simple.
  \end{thm}
  
  \parag{Proof} The equality $\mathcal{J} = [\Delta^{-1}] \mathcal{I} \cap W_2$ is clear thanks to Lemma \ref{14/8} below.\\
   Let $Q \in W_2$ such $Q(\delta) \not= 0$. There exists $m \in \mathbb{N}$ such that $\Delta^m Q$ is in $W_1^{\frak{S}_k}$ and satisfies $\Delta^mQ(\delta) \not= 0$. So $\Delta^mQ$ is not in $\mathcal{I}$ and then we have
  $$W_1^{\frak{S}_k} = \mathcal{I} +W_1^{\frak{S}_k}\Delta^mQ. $$
  So there exists $\Pi \in \mathcal{I}$ and $P \in W_1^{\frak{S}_k}$ such that 
  $$ 1 = \Pi + P\Delta^mQ .$$
  So $\mathcal{J} + W_2Q = W_2$ and the theorem is proved.$\hfill \blacksquare$\\

 \section{Antisymmetric PDO and symmetric PDO}
 
 \subsection{Symmetric and antisymmetric vector fields}
 
We recall here some elementary facts.
 
  \begin{lemma}\label{Sym VF}
  Any symmetric vector field in $W_1^{\frak{S}_k}$ is in the $\C[\sigma]$-module generated by $V_{p, 1}$ for $p \in [0, k-1]$.
  \end{lemma}
  
  For a proof see for instance \cite{[B.22]} Lemma 6.1.1.\\

 The anti-symmetric vector fields in $W_1$ are described by the following lemma.
  \begin{lemma}\label{anti-vect.}
For e	ach $h \in [1, k]$ $\delta \partial_h$  is a vector fiels in $W_1$ (anti-symmetric, of course). Moreover, any anti-symmetric vector field is of the form $\delta V$ where $V$ is a vector field in $W_2$. 
  \end{lemma}
  
  \parag{Proof} Using Theorem \ref{24/2/24} below we obtain that if $A$ an anti-symmetric vector field  then $V := \delta^{-1}A$ is a vector field in $W_2$. The converse is a consequence of the formula
    \begin{equation}
(-1)^{h-1}\partial_h = \sum_{j=1}^k \frac{z_j^{k-h}}{P'(z_j)}\frac{\partial}{\partial z_j} 
\end{equation}
which is proved in \cite{[B.22]} Lemma 6.1.2,  since $\delta$ is a multiple of $P'(z_j)$ in $\C[z_1, \dots, z_k]$ for each $j \in [1, k]$.$\hfill \blacksquare$\\

Let us complete this proof in giving the explicit formula for the anti-symmetric vector fields $\delta \partial_h$.\\
   Writing $\delta = (-1)^{j-1}P'(z_j)\vartheta(j) $ where 
  $$\vartheta(j) := \prod_{1 \leq i < h \leq k; i, j \not= j} (z_i - z_h)$$
   the above formula shows that
  $$ (-1)^{h-1}\delta\partial_h =  \sum_{j=1}^k  (-1)^{j-1}z_j^{k-h}\vartheta(j)\frac{\partial}{\partial z_j}  $$
  and the right hand-side is clearly in $W_1$ and then is an antisymmetric vector field in $W_1$.$\hfill \blacksquare$\\
  
  An easy consequence of the previous result is the following.
  
\begin{lemma}\label{14/8}
Let $P \in W_2$ be of order $q \geq 1$. Then $\delta^{2q-1}P$ is in $W_1^{\frak{A}_k, -}$ the subspace of   anti-symmetric differential operators in $W_1$.
\end{lemma}

An obvious consequence is that $\Delta^q P$ is in $W_1^{\frak{S}_k}$ if $P \in W_2$ has order $q$.

\parag{Proof} When $q = 1$ this consequence of Lemma  \ref{anti-vect.}. So assume the lemma proved for $q-1$ with $q \geq 2$ and consider $P$ of order $q$. Then write
$$ P = \sum_{h=1}^k \partial_h Q_h  + Q_0$$
where each $Q_h$ has order $\leq q-1$. Then since we have for each $h \in [1, k]$:
$$ \delta^{2q-1}\partial_h = \partial_h\delta^{2q-1} - (2q-1)\delta^{2q-2} \partial_h(\delta)$$
and so
$$ \delta^{2q-1}\partial_h Q_h = \partial_h\Delta \delta^{2q-3}Q_h - (2q-1)\delta\partial_h(\delta) \delta^{2q-3}Q_h .$$
 Our induction hypothesis implies that  $\delta^{2q-3}Q_h$ and  $\delta^{2q-1}Q_0 = \Delta \delta^{2q-3}Q_0 $ are in $W_1^{\frak{A}_k, -}$, so the conclusion follows using the fact that $2\delta\partial_h(\delta) = \partial_h(\Delta) $ is in $\C[\sigma]$ and  that $\partial_h\Delta = \Delta\partial_h + \partial_h(\Delta)$ is in $W_1^{\frak{S}_k}$.$\hfill \blacksquare$\\

 \subsection{Anti-Symmetric PDO and discriminant}
 
 The goal of this subsection is to investigate  the image  $\mathcal{M}_0 $ of $W_2/\mathcal{J}$ inside the $W_2$-module $\C[\sigma][\Delta^{-1}]\delta$ associated to the (regular)  connection defined by
 $$\nabla_{\partial_j}\delta = (1/2)\Delta^{-1}\partial_j(\Delta)\delta .$$
 We obtain, for instance, the fact  that $(1/\Delta)\delta$ belongs to $\mathcal{M}_0 \simeq W_2\delta$.\\
 
 A key tool is the following result.
 
 \begin{thm}\label{24/2/24}
 Let $A$ be an anti-symmetric differential operator in $W_1$. Then $\delta A$ is in $\Delta W_2$, where $\Delta = \delta^2$ is the discriminant. So any such $A$ is in $\delta^{-1}W_2$.
 \end{thm}
 
 \parag{Proof} Note first that it is enough to show this result for $k = 2$ because at the generic point of $\{ \Delta = 0\}$, the ramification set of the quotient map $q : M \to N$ by the action of $\frak{S}_k$ on $M = \C^k$ we have a decomposition of $q$ as the product of the quotient map $\C^2 \to \C^2/\frak{S}_2$ with an étale covering.\\
 Moreover, if we may write $\delta A = \Delta Q$ with $Q $ in $W_2$ localized near a point in $\{\Delta = 0 \} $, then $Q$ is unique, and since the sheaf associated to $W_2$ on $ \C^k$ satisfies the analytic extension property in co-dimension $\geq 2$ as it is an increasing union of free finite type $\mathcal{O}_{\C^k}$-modules, it is enough to show our result near the generic point of $\{\Delta = 0 \}$. So it is enough to prove the theorem for $k = 2$.\\
 In the case $k = 2$ note $z_1 = a$ and $z_2 = b$ and consider a monomial $a^pb^q(\partial_a)^r(\partial_b)^s$. Then $a^pb^q(\partial_a)^r(\partial_b)^s - a^qb^p(\partial_a)^s(\partial_b)^r$ is an anti-symmetric differential operator and clearly any anti-symmetric differential operator is a finite sum of such operators. \\
 So it is enough to prove the result for these special cases. \\
 But writing for $p \leq q$ and $r \leq s$ or $r > s$  without lost of generilty we have
 \begin{align*}
 & a^pb^q(\partial_a)^r(\partial_b)^s  - a^qb^p(\partial_a)^s(\partial_b)^r= \sigma_2^p\big(b^{q-p}(\partial_b)^{s-r} - a^{q-p}(\partial_a)^{s-r}\big)\Sigma_2^r \quad {\rm or}\\
 & a^pb^q(\partial_a)^r(\partial_b)^s - a^qb^p(\partial_a)^s(\partial_b)^r = \sigma_2^p\big(b^{q-p}(\partial_a)^{r-s} - a^{q-p}(\partial_b)^{r-s}\big)\Sigma_2^s
 \end{align*} 
  where we denote $\sigma_1 := a + b$, $\sigma_2 := ab$, $\Sigma_1 := \partial_a + \partial_b$ and $\Sigma_2 := \partial_a\partial_b$. \\
  Denote $\partial_1$ and $\partial_2$ the partial derivative in the  coordinates $\sigma_1, \sigma_2$ of $\C^2\big/\frak{S}_2$. Using the fact  that $\Sigma_1 = 2\partial_1 - \sigma_1\partial_2$ and  $\Sigma_2 = \partial_1^2 + 2\sigma_1\partial_1\partial_2 + 2\sigma_2\partial_2^2 + \partial_2$, it is enough to consider the anti-symmetrizations of the monomials 
 $$ a^p(\partial_a)^q \quad {\rm and} \quad a^p(\partial_b)^q \quad {\rm for \ all} \quad p, q \in \mathbb{N}.$$
Then we have
 \begin{align*}
 & a^p = x_pa + y_p \quad {\rm for} \ p \geq 2 \quad {\rm with} \ x_2 = \sigma_1 \quad {\rm and} \ y_2 = -\sigma_2 \tag{1} \\
 & {\rm where} \quad x_{p+1} = x_p\sigma_1 + y_p \quad {\rm and} \quad y_{p+1} = -x_p\sigma_2 
 \end{align*}
 and the analog formulas for $b^p$. But we have also
 \begin{align*}
 & \partial_a^p = X_p\partial_a + Y_p \quad {\rm with} \ X_2 = \Sigma_1 \quad {\rm and} \ Y_2 = -\Sigma_2 \tag{2} \\
  & {\rm where} \quad X_{p+1} = X_p\Sigma_1 + Y_p \quad {\rm and} \quad Y_{p+1} = -X_p\sigma_2 
  \end{align*}
   and the analog formulas for $\partial_b^p$. Note that $X_p$ and $Y_p$ commute with $\partial_a$ and $\partial_b$ and that they are in the commutative algebra generated by $\Sigma_1$ and $\Sigma_2$ (whose  elements commute with $\partial_a$ and $\partial_b$).\\
   Then we have the following special cases of our theorem
   \begin{align*}
   & (a-b)(\partial_a - \partial_b) = \Delta \partial_2 \\
   & (a-b)(a\partial_a -b\partial_b) = \Delta \partial_1 \\
   & (a-b)(b\partial_a -a\partial_b) = -\Delta(\partial_1 - \sigma_1\partial_2)
   \end{align*}
   where $\Delta = \delta^2 = (a - b)^2 = \sigma_1^2 - 4\sigma_2$.\\
     These cases correspond to the anti-symmetrizations of the monomials $ a^p(\partial_a)^q $ and $a^p(\partial_b)^q $ respectively for the cases $p = 0, 1$ and $q=0$ for the
    first one  and $p=1, q=1$ for the second one (note the case $p=0$ for the second one  is the same than  $p=0$ for the first one up to a sign).\\
   So consider now first the cases $p \geq 2$ for $a^p\partial_a$ or for  $a^p\partial_b$. The relation $(1)$ will allow us to reduce these case to $p =0$ and $p = 1$.\\
   In the same way the cases $q\geq 2$ for $a(\partial_a)^q$ or for $a(\partial_b)^q$ the relation $(2)$ will  allow us  to reduce these case to $q =0$ and $q = 1$.\\
   Assuming now that $p \geq 2$ and $q\geq 2$ we have
   \begin{align*}
   & a^p(\partial_a)^q - b^p(\partial_b)^q = (x_pa + y_p)\big(\partial_aX_q + Y_q\big) - (x_pb + y_p)\big(\partial_bX_q + Y_q\big) \\
   & \quad = x_p(a\partial_a  -b\partial_b)X_q + x_p(a-b)Y_q + y_p\big(\partial_a -\partial_b\big)X_q 
   \end{align*}
   and in analogous way
   \begin{align*}
   & a^p(\partial_b)^q - b^p(\partial_a)^q = (x_pa + y_p)\big(\partial_bX_q + Y_q\big) - (x_pb + y_p)\big(\partial_aX_q + Y_q\big) \\
    & \quad = x_p(a\partial_b  -b\partial_a)X_q + x_p(a-b)Y_q + y_p\big(\partial_b -\partial_a\big)X_q 
    \end{align*}
    So, after product by $\delta = (a-b)$ we see from the cases above that we find element in $\Delta W_2$ where here $W_2$ is the Weyl algebra $\C<\sigma_1, \sigma_2, \partial_1, \partial_2>$,
    concluding the proof of the theorem.$\hfill \blacksquare$\\
    
    \parag{Remark}  For any $P \in W_1^{\frak{S}_k}$ then $P\delta$ is anti-symmetric, then the previous result gives $Q \in W_2$ such that $\delta P\delta = \Delta Q$ so $\delta^{-1}P\delta = Q$. This shows that $ \delta^{-1}W_1^{\frak{S}_k}\delta \subset W_2$. Note that $ \delta^{-1}W_1^{\frak{S}_k}\delta$ is clearly a sub-algebra of $[\Delta^{-1}]W_2$.\\
   
It is not true in general  that an element $P \in W_1^{\frak{S}_k}$ which is in $\Delta W_2$ has its   coefficients (as an element in $W_1$)  which vanish on  $\{\delta = 0\}$. Let us give an example.
 
 \parag{Example} We consider the case $k = 2$ and we keep the previous notations. using the equality $\partial_2 =  (1/(a-b))(\partial_a - \partial_b)$ we obtain
 \begin{align*}
 & (a-b)^2\partial_2^2 = \frac{2}{a-b}(\partial_a -\partial_b) + (\partial_a - \partial_b)^2 \\
 & (a-b)^4\partial_2^3 = \frac{12}{a-b}(\partial_a -\partial_b) + 6((\partial_a - \partial_b)^2 + (a-b)(\partial_a - \partial_b)^3 
 \end{align*}
 and so $6\Delta\partial_2^2 - \Delta^2\partial_2^3 = -(a-b)(\partial_a - \partial_b)^3$.\\
  Now use the equality $\partial_1 = (1/(a-b))(a\partial_a - b\partial_b)$ we obtain
 $$ \Delta\big(6\partial_2^2 - \Delta\partial_2^3\big)\partial_1 = -(\partial_a - \partial_b)^3(a\partial_a - b\partial_b) .$$ 
 This gives an example of a $P \in W_1^{\frak{S}_k}$   which satisfies $P = \Delta Q $ with $Q \in W_2$ and such that the   coefficients of $P$ (as an element in $W_1$) are not vanishing on $\{\delta = 0\}$: \\
   the coefficient of $\partial_a^4$ in the right hand-side above is equal to  $-a$ which does not vanish identically  when $a = b$. So $\delta^{-1}P$ is anti-symmetric but not in $W_1$$\hfill\square$\\
 
 So it is not true that the image in $\Delta W_2$ of the symmetric differential operators of the form $\delta A$ with $A$ anti-symmetric in $W_1$  is equal to $\Delta W_2$.\\
  But this is true for vector fields (see above).
 
 \subsection{The computation of $\check{\delta}(\delta)$}
  Note $\delta_k := \prod_{1 \leq p < q\leq k} (z_p - z_q) $  the  discriminant of $z_1, \dots, z_k$ and  $\Cech{\delta}_k$  the  discriminant of $\partial_{z_1}, \dots, \partial_{z_k}$.\\
 We begin by some easy lemmas.

\begin{lemma}\label{0}
Let $P \in \C[\sigma]$ with weight $h \in [0, k]$ such that $P$ has degree  at most $1$ in $z_k$. Then $P = \alpha \sigma_h$ for some $\alpha \in \C$.
\end{lemma}

\parag{Proof} Write $P =  Q(z')z_k + R(z')$ where $Q$ and $R$ are $\frak{S}_{k-1}$-invariant of degrees $h$ and $h-1$ respectively. Now, by the $\frak{S}_k$-invariance of $P$ they are of degree $1$ at most in $z_{k-1}$. So assuming that the result is proved for $k-1$ we obtain for $h \leq k-1$  the equality  $P = \alpha \sigma_{h-1}(z')z_k + \beta \sigma_{h}(z')$ for some complex numbers $\alpha$ and $\beta$. The $\frak{S}_k$-invariance of $P$ implies then that $\alpha = \beta$ and the conclusion follows since $\sigma_h(z) = \sigma_h(z') + \sigma_{h-1}(z')z_k$.\\
For $h = k$ the $\frak{S}_{k-1}$-invariant polynomial  $R(z')$ has weight $k$ and degree at most $1$ in each variable $z_1, \dots, z_{k-1}$. So it must vanish and we have
$ P = Q(z')z_k$ where the $\frak{S}_{k-1}$-invariant polynomial $Q$ has weight $k-1$ and degree at most $1$ in each variable $z_1, \dots, z_{k-1}$. So $Q(z') = \alpha z_1\cdots z_{k-1}$ and we conclude that $P = \alpha \sigma_k$.$\hfill \blacksquare$\\

\begin{cor}\label{1}
Let $P$ in $\C[\sigma]$ of weight $h \in [1, k]$ such that $P\delta_k$ has degree $k$ in $z_k$, then $P = \alpha \sigma_h$.
\end{cor}

\parag{Proof} Since $\delta_k$ is a polynomial of degree $k-1$ in $z_k$, the previous lemma applies to $P$ which has degree $\leq 1$ in $z_k$.$\hfill \blacksquare$\\

\begin{lemma}\label{2}
For each $h \in [1, k]$ we have
\begin{equation*}
\Sigma_h[\sigma_k\delta_k] = h!\sigma_{k-h}\delta_k.\tag{F1}
\end{equation*}
\end{lemma}

\parag{Proof} Since $\Sigma_h[\sigma_k\delta_k] $ is anti-symmetric with weight $k-h + k(k-1)/2$ \  it can be written as $P\delta_k$ with $P$ of weight $k-h$. But the degree in $z_k$ of $\sigma_k\delta_k$
is equal to $k$ and $\Sigma_h$ derived at most one time in $z_k$, so the degree in $z_k$ of $P$ is at most $1$. This proves, thanks to the previous corollary,  the existence of a constant $\gamma(h,k)$ such that
$$ \Sigma_h[\sigma_k\delta_k] = \gamma(h,k)\sigma_{k-h}\delta_k $$
holds true for any $k \geq 1$ and any $h \in [1, k]$.\\
We shall use the equality:
\begin{equation*}
 \delta_k(z) = (-1)^{k-1}\Pi_k(z_k)\delta_{k-1}(z') \quad {\rm where} \quad \Pi_k(z) := \prod_{j=1}^{k-1} (z - z_j)  \tag{E}
 \end{equation*}
First for $h = k$ we see that the degree $k-1$ term  in $z_k$ inside  $\Sigma_k[\sigma_k\delta_k]$ is given by
$$ k \Sigma_{k-1}(z')[(-1)^{k-1}\sigma_{k-1}(z')\delta_{k-1}(z')] $$
since $\sigma_k\delta_k$ is  a degree $k$  in $z_k$ with leading  coefficient  equal to $(-1)^{k-1}\sigma_{k-1}(z')\delta_{k-1}(z')$ and, since  $\Sigma_k = \Sigma_{k-1}(z')\partial/\partial z_k $,  we find that
$ \gamma(k, k) = k\gamma(k-1, k-1)$ and so
$\Sigma_k[\sigma_k\delta_k] = k! \delta_k$ since $\gamma(1, 1) = 1$.\\
But for $h \leq k-1$  the degree $k$ term  in $z_k$ inside  $\Sigma_h[\sigma_k\delta_k]$  is given by
$$ \Sigma_h(z')[ (-1)^{k-1}\sigma_{k-1}(z')\delta_{k-1}(z')] .$$
 using the equality $(E)$.

This gives the relation, since the coefficient of $z_k^k$ in $\sigma_k(z)\delta_k(z)$ is equal to $\sigma_{k-1}(z')(-1)^{k-1}\delta_{k-1}(z')$:
$$  \Sigma_h(z')[ (-1)^{k-1}\sigma_{k-1}(z')\delta_{k-1}(z')]  = \gamma(h,k) \sigma_{k-1}(z')(-1)^{k-1}\delta_{k-1}(z')$$
and then $\gamma(h, k-1) = \gamma(h, k)$.\\
 The conclusion follows  from the fact that  $\gamma(k,k) = k!$.$\hfill \blacksquare$\\
 
\begin{lemma}\label{3}
For any  $h \in [1, k]$ and any $p \in [1, h]$ we have
\begin{equation*}
 \Sigma_1^p[\sigma_h\delta_k] = \frac{(k-h+p)!}{(k-h)!}\sigma_{h-p}\delta_k. \tag{F2}
 \end{equation*}
\end{lemma}

\parag{Proof} This is an obvious consequence of the equalities  $\Sigma_1[\delta_k] = 0$ and  \\
 $\Sigma_1[\sigma_q] = (k-q+1)\sigma_{q-1}$ using the Leibnitz rule for a vector field.$\hfill \blacksquare$\\

\begin{prop}\label{enfin !}
For each $ k \geq 2$ we have the formula
\begin{equation*}
\Cech{\delta}_k[\delta_k] = c_k \tag{F0}
\end{equation*}
where $c_k$ is a positive constant equal to $c_k = k!(k-1)!\dots 2!$.
\end{prop}

\parag{Proof} It is clear, looking to the weight of the left hand-side in $F(k)$ that the result is a constant, since it is a weight $0$ polynomial. So we are looking for the constant $c(k)$ of the right hand-side.\\
First for $k = 2$ we have $(\partial_a - \partial_b)[a - b] = 2(a - b)$ so $c_2 = 2$.\\
Now we shall argue by induction on $k \geq 2$. Looking a	t the change of variable given by
$$ x_j = z_j -z_{k+1} \quad {\rm for} \quad j \in [1, k] \quad {\rm and} \quad x_{k+1} = z_{k+1} $$
we obtain that 
$$(-1)^k \delta_{k+1}(z) = x_1\dots x_k \delta_k(x) = \sigma_k(x) \delta_k(x)$$
and 
$$ \Cech{\delta}_{k+1}(z) = (-1)^k\Cech{\delta}_k(x)\prod_{j=1}^k\big(\partial_{x_j} + \Sigma_1(x) - \partial_{x_{k+1}}\big).$$
So this gives, since we may  omit $\partial_{x_{k+1}}$ because $\big[x_1\dots x_k \delta_k(x) \big]$ does not depend on the variable $x_{k+1}$
\begin{align*}
& \Cech{\delta}_{k+1}(z)[\delta_{k+1}(z)] = \Cech{\delta}_k(x)\prod_{j=1}^k(\partial_{x_j} + \Sigma_1(x))\big[x_1\dots x_k \delta_k(x) \big] \\
& c_{k+1} = \Cech{\delta}_k(x)\prod_{j=1}^k(\partial_{x_j} + \Sigma_1(x))\big[x_1\dots x_k \delta_k(x) \big]
\end{align*}
Now  we have 
$$ \prod_{j=1}^k (X + \xi_j) = \sum_{h=0}^k \sigma_h (\xi)X^{k-h}$$
and this  implies
$$ \prod_{j=1}^k(\partial_{x_j} + \Sigma_1(x)) = \sum_{h=0}^k \Sigma_h(x)(\Sigma_1(x))^{k-h} =   \sum_{h=0}^k (\Sigma_1(x))^{k-h}\Sigma_h(x)$$
Now using the formula $(F1)$ and $(F2)$ we find
\begin{align*}
& c_{k+1} = \Cech{\delta}_k(x)\Big[ \sum_{h=0}^{k}  (\Sigma_1(x))^{k-h}\big[\Sigma_h(x)[\sigma_k(x)\delta_k(x)]\big] \Big]\\
& \quad = \Cech{\delta}_k(x)\Big[ \sum_{h=0}^{k}  (\Sigma_1(x))^{k-h}\big[h! \sigma_{k-h}(x)\delta_k(x)\big] \Big]\\
& \quad =  \big( \sum_{h=0}^{k} k! \big) c_k =  (k+1)!  c_k
\end{align*}
Since $c_2 = 2$ the proof is complete.$\hfill \blacksquare$\\

\begin{cor} \label{quid}
There exists $P \in W_2$ such that $P(\delta) = 1/\delta = (1/\Delta)\delta$.
\end{cor}

\parag{Proof} Since $\delta\Check{\delta} = \Delta P_1$ for some $P_1 \in W_2$, thanks to Theorem \ref{24/2/24}, we obtain $\Delta P_1(\delta) = c_k\delta$ so $P := c_k^{-1}P_1$ satisfies the relation  $P(\delta) = 1/\delta = (1/\Delta)\delta$.$\hfill \blacksquare$\\

\parag{Remark} The fact that $\check{\delta}(\delta) = c_k \not= 0$ shows that $\check{\delta}$ does not belong to $\delta W_1^{\frak{S}_k}$ since $\check{\delta}$ does not send $\C[\sigma]\delta$ to $\C[\sigma]\delta$ !\\
Note that $\delta \not\in W_1^{\frak{S}_k}\check{\delta}$ since $\check{\delta}(1) = 0$.

   \begin{cor}\label{29/9/24}
Let $b$ be the Bernstein polynomial of the discriminant $\Delta := \delta_k^2$ in $\C[\sigma]$. Then $b(-1/2) \not= 0$.
\end{cor}

\parag{Proof}This is a simple consequence of the existence of $Q \in W_2$ such that $Q(\Delta^{1/2}) = \Delta^{-1/2}$ and the following result which is a simple consequence of the definition of the Bernstein polynomial (see \cite{[K.76]}).

 \begin{prop}\label{K.76}
Let $(f, 0) : (\C^n, 0) \to (\C, 0)$ be a non zero  germ of holomorphic function and let $X$ be a small open neighborhood of $0$ on which $f$ is defined and satisfies
 $\{df = 0\} \subset \{f = 0\}$.\\
  Let $\mathcal{N}_f := D_X[s]f^s$ be the $D_X[s]$-module generated by $f^s$ inside $\mathcal{O}_X[s, f^{-1}]f^s$ and  define $$t \in \mathcal{H}om_{D_X}(\mathcal{N}_f, \mathcal{N}_f) \quad {\rm by} \quad  t(P(s)f^s = P(s+1)f^{s+1} = P(s+1)f.f^s.$$
  The Bernstein polynomial $b_{f, 0}$  of $f$ at the origin is, by definition,  the minimal polynomial of the germ of \  $t$  \ at the origin and then $b_{f, 0}(r) \not= 0$ for some complex number $r$ if and only if the germ of \ $t$  \ at $0$ induces an isomorphism
 of $\mathcal{N}_f/(s-r)\mathcal{N}_f \simeq D_X f^r$ onto itself. So $b_{f, 0}(r) \not= 0$  is equivalent to the existence of a germ $Q \in D_{X, 0}$ which satisfies $Qf^{r+1} = f^r$.$\hfill \blacksquare$\\
 \end{prop}


\bigskip
\bigskip

  \section{Bibliography}
   
    \begin{thebibliography}{99}
    
    \bibitem{[K.76]} Kashiwara, M. {\it B-Functions and Holonomic Systems},  Inventiones  math. 38, (1976) pp. 33-53.
    \bibitem{[B.21]} Barlet, D. {\it Note on Lisbon integrals and their associated  D-modules}, Port. Math. (N.S.) 78 $n^0$ 3-4 (2021), pp. 323-340.
    \bibitem{[B.22]} Barlet, D. {\it On symmetric partial differential operators}, Math. Z. 302, $n^0 3$ (2022) pp. 1627-1655.


\end{thebibliography}
\end{document}